\input amstex
\documentstyle{amsppt}
\magnification 1200

\NoRunningHeads
\topmatter
\title Martingale-coboundary decomposition for stationary  random fields
\endtitle
\author
Dalibor Voln\'y
\endauthor
\abstract
We prove a martingale-coboundary representation for random fields with a completely commuting filtration. 
For random variables in $L^2$ we present a necessary and sufficient condition which is a generalization of Heyde's condition
for one dimensional processes from 1975. For $L^p$ spaces with $2\leq p<\infty$ we give a necessary and sufficient condition which 
extends Voln\'y's result from 1993 to random fields and improves condition of El Machkouri and Giraudo from 2016. A new sufficient condition is presented
which for dimension one improves Gordin's condition from 1969. In application, new weak invariance principle and estimates of large deviations are found.
\endabstract
\endtopmatter

\document 
\subheading{1. Introduction}

Let $(\Omega, \Cal A, \mu)$ be a probability space and $(T_{\underline{i}})_{\underline{i}\in \Bbb Z^d}$ a $\Bbb Z^d$ action on $(\Omega, \Cal A, \mu)$
generated by  commuting invertible and measure-preserving transformations $T_{\epsilon_q}$, $1\leq q\leq d$.
By $\epsilon_q$ we denote the vector from $\Bbb Z^d$ which has 1 at $q$-th place and 0 elsewhere. 
By $U_{\underline{i}}$ we denote the operator in $L^p$ ($1\leq p<\infty$) defined by $U_{\underline{i}}f = f\circ T_{\underline{i}}$,
$\underline{i}\in \Bbb Z^d$. 

By $\underline{i} \leq \underline{j}$ we understand $i_q\leq j_q$ for all $1\leq q \leq d$.
The vectors $(0,\dots,0)$ and $(1,\dots,1)$ will be denoted $\underline{0}$ and $\underline{1}$ respectively.

We suppose that there is a {\it completely commuting filtration} $(\Cal F_{\underline{j}})_{\underline{j}\in\Bbb Z^d}$, i.e\. there is a $\sigma$-algebra 
$\Cal F$ such that $\Cal F_{\underline{i}} = T_{-\underline{i}}\Cal F$, for  $\underline{i}\leq \underline{j}$ we have  
$\Cal F_{\underline{i}} \subset \Cal F_{\underline{j}}$, and for an integrable $f$ it is
$$
  E\big(E(f\, |\, \Cal F_{i_1,i_2,\dots,i_d}) \,|\, \Cal F_{j_1,j_2,\dots,j_d}\big) = 
  E(f\,|\, \Cal F_{i_1\wedge j_1,i_2\wedge j_2,\dots,i_d\wedge j_d})
$$
where $i\wedge j = \min \{i, j\}$ (cf\. \cite{VWa14}). As a frequent (cf\. \cite{VWa14}, \cite{WaWo13} and references therein) but not exclusive 
example, let us introduce  a Bernoulli $\Bbb Z^d$ action:  the $\sigma$-algebra $\Cal A$ is 
generated by iid random variables $e_{\underline{i}} = U_{\underline{i}} e$. The filtration $\Cal F_{\underline{j}} = \sigma\{e_{\underline{i}} : 
\underline{i} \leq \underline{j}\}$ is completely commuting.

By $\Cal F_l^{(q)}$ we denote the $\sigma$-algebra generated by all $\Cal F_{\underline{i}}$ with $i_q\leq l$ ($i_j\in\Bbb Z$
for $1\leq j\leq d$, $j\neq q$), $1\leq q\leq d$. \newline
For $\sigma$-algebras $\Cal G\subset \Cal F \subset \Cal A$ and $1\leq p<\infty$ we by $L^p(\Cal F)\ominus L^p(\Cal G)$ denote the space 
of $f\in L^p(\Cal F)$ for which $E(f\,|\,\Cal G)=0$.
Similarly as in the one dimensional case we can define projection operators $P_l^{(q)}$ onto 
$L^p(\Cal F_l^{(q)})\ominus L^p(\Cal F_{l-1}^{(q)})$ by 
$P_l^{(q)}f = E(f\,|\,\Cal F_l^{(q)}) -  E(f\,|\,\Cal F_{l-1}^{(q)})$. These operators commute and for $l\neq k$,  $P_l^{(q)} P_k^{(q)} = 0$.
We define  projections $P_{j_1,j_2,\dots,j_d} =  P_{j_1}^{(1)} \dots  P_{j_d}^{(d)}$ onto 
$\underset 1\leq q\leq d \to{\bigcap} L^p(\Cal F_{j_1}^{(q)}) \ominus  L^p(\Cal F_{j_1-1}^{(q)})$ (cf\. \cite{VWa14}).
Let us notice that $T_{\epsilon_q}^{-1} \Cal F_{j_1,j_2,\dots,j_d} = \Cal F_{j'_1,j'_2,\dots,j'_d}$ where $j_q' = j_q+1$ and $j'_i=j_i$ for $i\neq q$.
We thus have $T_{\epsilon_q}^{-1} \Cal F_l^{(q)} =
\Cal F_{l+1}^{(q)}$ and  $T_{\epsilon_q}^{-1} \Cal F_l^{(q')} = \Cal F_l^{(q')}$ for $q'\neq q$; $U_{\epsilon_q} P_l^{(q)} = 
P_{l+1}^{(q)}U_{\epsilon_q}$ and $U_{\epsilon_q} P_l^{(q')} =  P_l^{(q')} U_{\epsilon_q}$. 

 An integrable function $f$  is called {\it regular}  if it is $\Cal{ F}_{\underline{\infty}}$-measurable and for every 
$1\leq q\leq d$, 
$E(f\,|\, \Cal{ F}_{-\infty}^{(q)}) = 0$, where $\Cal{ F}_{-\infty}^{(q)} = \cap_{k\in\Bbb Z} \Cal{ F}_{-k}^{(q)}$.
A function $f$ is {\it adapted} if it is $\Cal F_{\underline{0}}$-measurable.
$f\in L^p$, $1\leq p<\infty$, is thus regular if and only if $f= \sum_{\underline{i}\in\Bbb Z^d}  P_{\underline{i}}f$.
In this paper all functions will be supposed to be regular. A random field $U_{\underline{i}}f$ generated by a regular function $f$ will be called
regular.

A useful tool in proving limit theorems for one dimensional (strictly) stationary processes $(f\circ T^i)_i$  (i.e\. for $d=1$) has been the 
martingale-coboundary decomposition
$$
  f = m + g - g\circ T \tag{1}
$$
where $(m\circ T^i)_i$ is a martingale difference sequence. 
The decomposition (1) was one of the first conditions giving a CLT for stationary sequences of random variables by martingale approximation;
for $f, m, g\in L^2$ it was introduced already in Gordin's 1969 paper \cite{Go69}.
Even if the cobounding function $g$ is just measurable, (1) still implies a central limit theorem. A paper proving a CLT for $f,g\in L^1$ was proved
in \cite{Go73} (cf\. also \cite{EJ85}; see that square integrability of the martingale differences $m$ needed to be proved).
Square integrability of $m$, $g$ guarantees a weak invariance principle (WIP) and a functional law of iterated logarithm (\cite{He75}).
Notice that in general, a central limit theorem does not imply WIP. For strictly stationary and ergodic processes this has been shown e.g\. in \cite{VSa00}.
 In \cite{GV14} a beta mixing process satisfying the CLT but not WIP is found. 

The condition (1) provides a very close martingale approximation and for central limit theorems it is sometimes suboptimal, e.g\. the conditions of 
Dedecker-Rio and of  Maxwell-Woodroofe imply the weak invariance principle (cf\. \cite{DR00}, \cite{MW00}).
The conditions mentioned above follow from (1) (with $m, g\in L^2$) but not vice versa, cf\.  e.g\. \cite{DuV08}.
(1) is independent of the Hannan's criterium (cf\. \cite{V93}) and remains useful in the study of central limit theorems for Markov chains (one of the first papers on the subject is \cite{GoLi}).
The martingale-coboundary decomposition (1) can be used in proving other limit theorems like estimates of large deviations
(cf\. \cite{LV01}) where other conditions do not apply. This motivates study of (1) in $L^p$ spaces with $p>2$. Probably the most 
exhaustive study of (1) in various spaces is in  \cite{V06}.

\medskip

In this paper we will extend the martingale-coboundary decomposition to random fields. 
In dimension $d\geq 2$ the decomposition appears more complicated:
we are interested in the existence of the (martingale-coboundary) representation
$$
  f = \sum_{S\subset \{1,\dots,d\}}\prod_{q\in S^c} (I-U_{\epsilon_q})g_S\,\,\,\,(f\in L^p, \,\, 1\leq p<\infty) \tag2
$$
where for $S\subset \{1,\dots,d\}$, $g_S \in \underset q\in S \to{\bigcap}  L^p(\Cal F_{0}^{(q)}) \ominus  L^p(\Cal F_{-1}^{(q)})$;
$\prod_{q\in \emptyset} (I-U_{\epsilon_q})$ is defined as $I$, the identity operator. For $q'\in S$, 
$U_{\epsilon_{q'}}^i\prod_{q\in S^c} (I-U_{\epsilon_q})g_S$, $i\in \Bbb Z$, are thus martingale differences while for $q'\in S^c$, 
$\prod_{q\in S^c} (I-U_{\epsilon_q})g_S$ are coboundaries for the transformation $T_{\epsilon_{q'}}$.

As an illustration, consider $d=2$. (2) then becomes 
$$
  f = m + \big[ g_1 - U_{1,0}g_1 \big] + \big[ g_2 - U_{0,1}g_2\big]  + \big[ g - U_{1,0}g - U_{0,1} (g - U_{1,0}g)\big] 
$$
where $m, g \in L^p$,  $P_{0,0} m=m$ (i.e\. $U_{i,j}m$ are martingale differences), 
$g_1 \in L^p(\Cal{ F}_{0}^{(2)}) \ominus  L^p(\Cal{ F}_{-1}^{(2)})$,
$g_2 \in L^p(\Cal{ F}_{0}^{(1)}) \ominus  L^p(\Cal{ F}_{-1}^{(1)})$. The term $g_1 - U_{1,0}g_1 \in L^p(\Cal{ F}_{0}^{(2)}) \ominus  
L^p(\Cal{ F}_{-1}^{(2)})$ 
is thus a martingale difference sequence for the transformation $T_{0,1}$ and a coboundary for $T_{1,0}$; the last term is a coboundary.
\medskip

The aim of this paper is to study both sufficient, and necessary and sufficient conditions for  the decomposition (2) for regular functions.
In the same setting as here, a sufficient condition was recently found by El Machkouri and Giraudo in \cite{ElG16}.
Their main result (cf\. Theorem 5 here) is a multiparameter version of Gordin's one-dimensional sufficient condition from 1969 (cf\. \cite{Go69}).
In 2009, the problem was studied by Gordin in \cite{Go09} where, instead of $\Bbb Z^d$, he used semigroup  $\Bbb Z_+^d$ and instead 
of martingale differences he got
 reversed martingale differences. For $d=1$ his condition becomes the Poisson equation. Our results can be easily
converted to the setting applied in Gordin's paper.

We prove a multiparameter version of necessary and sufficient conditions from \cite{He75} (Theorem 2) and \cite{V93} (Theorem 4).
In Theorem 6 we present a sufficient condition which seems to be easier to verify than assumptions of Theorem 2 and Theorem 4.

The martingale-coboundary representation will be used in proving limit theorems, in particular a weak invariance principle (WIP) and estimates of 
probabilities of large deviations. We will extend similar results from \cite{ElG16}.

\subheading{2. Main results}

Let us recall that in the paper we suppose regularity of the function $f$. We will need

\proclaim{Proposition 1} For $1\leq p<\infty$, the decomposition (1) with $f, m, g\in L^p$ is equivalent to the convergence
of 
$$
  \sum_{j= 0}^\infty E(U^jf | \Cal F_{-1}),\quad \sum_{j= 1}^\infty  \Big[ (U^{-j}f -  E(U^{-j}f | \Cal F_{-1}) \Big]  \tag3
$$  
in $L^p$. The transfer function $g$ can be regular and we can fix $m=P_0m$. In such a case 
$$
  m = \sum_{i\in\Bbb Z} P_0U^if, \quad g = \sum_{j= 0}^\infty E(U^jf | \Cal F_{-1}) - 
  \sum_{j= 1}^\infty  \Big[ (U^{-j}f -  E(U^{-j}f | \Cal F_{-1}) \Big].
$$
\endproclaim

A proof can be found in \cite{V93} (for $p=1,2$) and in \cite{V06}.

In the case of  $d\geq 2$ we will prove a necessary and sufficient condition for (2)  in $L^2$.

\proclaim{Theorem 2} The martingale-coboundary decomposition (2) holds in $L^2$  if and only if for every $S\subset \{1,\dots,d\}$ and
$S'\subset S^c$
$$
  \sum_{j_u\geq 1,\,u\in S'}  \sum_{j_v\geq 1,\,v\in S^c\setminus S'}
  \big\|   \sum_{i_r\in\Bbb Z,\, r\in S} \,\, \sum_{i_u\geq \,j_u,\,u\in S'} \,\,  \sum_{i_v\geq j_v,\,v\in S^c\setminus S'} 
   P_{\underline{0}} U_{i_r, i_u,-i_v} f  \big\|_2^2 < \infty ;  \tag{4a}
$$
this is equivalent to
$$
  \sum_{S\subset \{1,\dots,d\}} \sum_{j_k\geq 0,\,\,k\in S } \sum_{j_l\leq 0,\,\,l\in S^c} 
  \big\| \sum_{i_k\geq j_k,\,\,k\in S} \sum_{i_l\leq j_l,\,\,l\in S^c}  P_{\underline{0}}U_{i_1,\dots,i_d}f \big\|_2^2 
  < \infty.   \tag{4b}
$$
If (2), (4) are valid then the functions $g_S$ can be regular and we then get
$$\multline
  g_S = \\
  \sum_{S'\subset S^c} (-1)^{|S^c\setminus S'|}  \sum_{i_r\in\Bbb Z,\, r\in S}\,\, \sum_{j_u\geq 1,\,i_u\geq 0,\,u\in S'}  
  \sum_{j_v\geq 0,\,i_v\geq 1,\,v\in S^c\setminus S'} P_{0,-j_u,j_v}U_{i_r, i_u,-i_v} f.
\endmultline \tag{5}
$$
\endproclaim

In the formula (5), $P_{0,-j_u,j_v}$ is the projection operator $P_{\underline{j}}$ where $j_i=0$ for $i\in S$, $j_i\leq -1$ for $i\in S'$, and
$j_i \geq 0$ for $i\in S^c\setminus S'$; the expression (operator) $U_{i_r, i_u,-i_v}$ is to be understood similarly.

For $d=1$ the condition (4b) is equivalent to
$$
  \sum_{j=1}^\infty \| \sum_{i=j}^\infty  P_0 U^if \|_2^2 + \sum_{j=1}^\infty \| \sum_{i=j}^\infty  P_0 U^{-i}f \|_2^2 < \infty.\tag6
$$
The condition (6) was found as sufficient  for (1) by C.C\. Heyde (\cite{He75}, cf\. also \cite{HaHe, Theorem 5.5}) while in \cite{V93}  it was proved 
necessassary and sufficient.

From (5) and orthogonality it follows
$$\multline
  \|g_S\|_2^2 = \\
  \sum_{S'\subset S^c}  \sum_{j_u\geq 1,\,u\in S'}  \sum_{j_v\geq 0,\,v\in S^c\setminus S'}
  \Big\|  \sum_{i_r\in\Bbb Z} \sum_{i_u\geq 0,\,u\in S'} \sum_{i_v\geq 1,\,v\in S^c\setminus S'} P_{0,-j_u,j_v} U_{i_r, i_u,-i_v} f  \Big\|_2^2, 
  \endmultline 
$$
$S\subset \{1,\dots,d\}$.

Let us recall that we use regularity assumption; if we take all functions $g_S$ regular then they are unique.
\medskip

To prove Theorem 2 we can use a superlinear random field representation (cf\.  \cite{VWoZ11, Theorem 1} for $d=1$, \cite{CCo13} in the general case):
there exist  $e_k = P_{\underline{0}} e_k$ and real numbers 
$a_{k,\underline{i}}$ such that 
$$
  \|e_k\|_2 =1, \quad k\geq 0,\quad \sum_{k=0}^\infty \sum_{\underline{i}\in\Bbb Z^d} a_{k,\underline{i}}^2 < \infty,\quad\text{and}\quad
   f = \sum_{k=0}^\infty \sum_{\underline{i}\in\Bbb Z^d} a_{k,\underline{i}} U_{-\underline{i}}e_k. \tag{7}
$$ 
To see this, let's notice that without loss of generality we can suppose that the $\sigma$-algebra $\Cal A$ is countably generated. Then there exists a countable orthonormal
basis $\{e_k\,:\,k=0,1,\dots\}$ of  $\underset 1\leq q\leq d \to{\bigcap}  L^2(\Cal F_{0}^{(q)}) \ominus  L^2(\Cal F_{-1}^{(q)})$, and 
$U_{-\underline{i}}e_k$, $k\geq 0$, $\underline{i}\in \Bbb Z^d$, is the orthonormal basis of the Hilbert space of regular functions from $L^2$.

If $ f = \sum_{\underline{i}\in\Bbb Z^d} a_{\underline{i}} U_{-\underline{i}}e$, 
$e\in \underset 1\leq q\leq d \to{\bigcap}  L^2(\Cal F_{0}^{(q)}) \ominus  L^2(\Cal F_{-1}^{(q)})$, we will speak of a stationary linear field.

Let $L_k$, $k=0,1,\dots$, denote the Hilbert space generated by $U_{\underline{i}}e_k$, $\underline{i}\in \Bbb Z^d$, and $\Pi_k$ the orthogonal projection
operator on $L_k$. The space of regular elements from $L^2$ is the direct sum $\underset k\geq 0 \to{\oplus} L_k$. Each of the spaces $L_k$ is invariant
w.r.t\. $U_{\underline{i}}$, $\underline{i}\in \Bbb Z^d$, and $\Pi_k$ commutes with the operators   $U_{\underline{i}}$. In (2) we thus get
$$
   f = \sum_{k= 0}^\infty f_k =  \sum_{k\geq 0} \sum_{S\subset \{1,\dots,d\}}\prod_{q\in S^c} (I-U_{\epsilon_q})g_{k,S}
$$
where $f_k, g_{k,S} \in L_k$ and $g_S = \sum_{k=0}^\infty g_{k,S}$, $\|g_S\|_2^2 = \sum_{k=0}^\infty \|g_{k,S}\|_2^2$.

Heyde's condition (5) can be deduced from Proposition 1. We prove it in a form useful for proving Theorem 2.
\medskip

\proclaim{Proposition 3 (C.C\. Heyde)}  For $f\in L^2$ regular, (6) is equivalent to the martingale-coboundary decomposition (1) where
for $f$ represented by (7)
$$
  m= \sum_{k=0}^\infty \sum_{i\in \Bbb Z} a_{k,i} e_k, \quad  g= \sum_{k=0}^\infty  \sum_{j=1}^\infty  \sum_{i=j}^\infty a_{k,i} U^{-j}e_k 
  - \sum_{k=0}^\infty \sum_{j=0}^\infty  \sum_{i=j+1}^\infty a_{k,-i} U^{j}e_k.  
$$
\endproclaim

\demo{Proof} 

For $f$ represented by (7), (6) is equivalent to
$$
   \sum_{k=0}^\infty  \sum_{j=1}^\infty \big(\sum_{i=j}^{\infty} a_{k,i}\big)^2 +
  \sum_{k=0}^\infty  \sum_{j=1}^\infty \big(\sum_{i=j}^{\infty} a_{k,-i}\big)^2 <\infty. \tag{6b}
$$

For simplicity sake let's suppose that $f$ is adapted, i.e\.
$
  f = \sum_{k=0}^\infty \sum_{i\geq 0} a_{k,i}U^{-i}e.
$
By Proposition 1, (1) is then equivalent to the convergence  (in $L^2$) of 
$$
  \sum_{j= 0}^\infty E(U^jf | \Cal F_{-1}) = \sum_{i< 0} P_i \sum_{j= 0}^\infty E(U^jf | \Cal F_{-1}) = 
   \sum_{k=0}^\infty  \sum_{i=1}^\infty \Big( \sum_{j=0}^\infty a_{k,i+j}\Big) U^{-i}e_k , 
$$  
i.e\. to the convergence of 
$
  \sum_{k=0}^\infty  \sum_{i=1}^\infty \Big( \sum_{j=0}^\infty a_{k,i+j}\Big)^2,
$
which, for $f$ adapted, is equivalent to the condition (6b). \newline
The proof of the general non adapted case is similar.
\enddemo
\qed

\demo{Proof of Theorem 2}  We suppose that $f$ is represented by (7).

The condition (4a) is then equivalent to
$$
  \sum_{k=0}^\infty  \sum_{j_u> 0,\,\,u\in S' } \sum_{j_v<0,\,\,v\in S^c\setminus S'} 
   \big(\sum_{i_r\in \Bbb Z,\, r\in S} \sum_{i_u\geq j_u,\,\,k\in S'} \sum_{i_v\leq j_v,\,\,l\in S^c\setminus S}  a_{k,i_1,\dots,i_d} \big)^2  < \infty \tag{4c}
$$
and (4b) is equivalent to
$$
  \sum_{k=0}^\infty  \sum_{j_u\geq 0,\,\,u\in S } \sum_{j_v\leq 0,\,\,v\in S^c} 
   \big( \sum_{i_u\geq j_u,\,\,k\in S} \sum_{i_v\leq j_v,\,\,l\in S^c}  a_{k,i_1,\dots,i_d} \big)^2  < \infty. \tag{4d}
$$
By using elementary equality
$
  \sum_{j=0}^\infty \big(\sum_{i=j}^\infty a_i\big)^2 =  \big(\sum_{i=0}^\infty a_i\big)^2 + \sum_{j=1}^\infty \big(\sum_{i=j}^\infty a_i\big)^2 
$
and induction we can prove that the sum of all (4c)  over $S\subset \{1,\dots,d\}$ and $S'\subset S^c$ equals the sum of all (4d) over  
$S\subset \{1,\dots,d\}$. The conditions (4a) and (4b) are thus equivalent.

 (5) becomes
$$\multline
   g_S =\\
  \sum_{k=0}^\infty \sum_{i_r\in\Bbb Z} \sum_{S'\subset S^c} (-1)^{|S^c\setminus S'|} \sum_{j_u\geq 1,i_u\geq j_u, u\in S'} 
  \sum_{j_v\geq 0,i_v\geq j_v+1, v\in |S^c\setminus S'} a_{i_r,i_u,-i_v}U_{0,-j_u,j_v} e_k. 
\endmultline \tag{5a}
$$

Let us prove equivalence of (2) and (4).

For $d=1$, Theorem 2 becomes the Heyde's theorem mentioned above (cf\.  Proposition 3). 

Let's suppose that $d\geq 2$ and that for $d-1$ the theorem is true.

\medskip

Because the operators $U_{\underline{i}}$ commute with the projections $\Pi_k$ and $g_S  = 
\sum_{k=0}^\infty g_{k,S}$ ($= \sum_{k=0}^\infty \Pi_k g_S$), $\|g_S\|_2^2 = 
\sum_{k=0}^\infty \|g_{k,S}\|_2^2$, (2) holds in $L^2$ if and only if it holds in all spaces $L_k$.  It is thus
sufficient to prove the theorem for a stationary linear process.
\medskip

 For simplicity of notation we suppose that $f$ is adapted.
The expression (5a) then becomes
$$
  g_S = \sum_{i_r\geq 0, \,r\in S}  \sum_{j_u\geq 1,\,u\in S^c} \sum_{i_u\geq j_u\,u\in S^c}  a_{i_r,i_u} U_{0,-j_u} e, \tag{5b}
$$
$S\subset \{1,\dots,d\}$ and $(0, -j_u)$ is a vector from $\Bbb Z^d$; (4d) becomes 
$$
  \sum_{j_1=0}^\infty \dots \sum_{j_d=0}^\infty \Big(\sum_{i_1=j_1}^\infty \dots \sum_{i_d=j_d}^\infty a_{i_1,\dots,i_d} \Big)^2 <\infty. \tag{4e}
$$
\medskip

Let's suppose (4e) (we have a stationary adapted random field). We will prove (2). \newline
If we apply  Proposition 3 to $U_{\epsilon_d}$ and the filtration $(\Cal F_i^{(d)})_i$, we by using (4e) deduce
$$
  f = m_d + g_d - U_{\epsilon_d}g_d \tag{8}
$$
where
$$ \gather
  m_d = \sum_{i_1= 0}^\infty \dots \sum_{i_{d-1}= 0}^\infty  \sum_{i_{d}= 0}^\infty a_{i_1,\dots,i_{d-1},i_d} U_{-i_1,\dots-,i_{d-1},0}e, \\
  g_d = \sum_{i_1= 0}^\infty \dots \sum_{i_{d-1}= 0}^\infty  \sum_{j_{d}= 1}^\infty  \sum_{i_{d}= j_d}^\infty 
  a_{i_1,\dots,i_{d-1},i_d} U_{-i_1,\dots-,i_{d-1},-j_d}e.
  \endgather
$$
For $i_d\geq 0$ let us denote
$$
 f_{i_d} =  \sum_{i_1= 0}^\infty \dots \sum_{i_{d-1}= 0}^\infty a_{i_1,\dots,i_{d-1},i_d}  U_{-i_1,\dots-,i_{d-1},0}e.
$$
We thus have
$$
  m_d =  \sum_{i_d= 0}^\infty f_{i_d}, \quad g_d =  \sum_{j_d= 1} ^\infty  \sum_{i_d= j_d}^\infty U_{\epsilon_{d}}^{-j_d}  f_{i_d}.  
$$
When applying Theorem 2  to an action of $T_{\epsilon_1}, \dots, T_{\epsilon_{d-1}}$ (notice that (4e) remains satisfied) we get
$$ 
  f_{i_d} = \sum_{S\subset\{1,\dots,d-1\}}  \prod_{q\in S^c\setminus\{d\}} (I-U_{\epsilon_q})g_{S,i_d}
$$
where $g_{S,i_d}$ is defined by (5b) applied to an action of $T_{\epsilon_1}, \dots, T_{\epsilon_{d-1}}$. By (8) we thus have
$$\multline
  f = m_d + (I-U_{\epsilon_d})g_d = \sum_{S\subset\{1,\dots,d-1\}}  \prod_{q\in S^c\setminus\{d\}} (I-U_{\epsilon_q})  \sum_{i_d= 0}^\infty g_{S,i_d}+ \\
  \sum_{S\subset\{1,\dots,d-1\}}   (I-U_{\epsilon_d}) \prod_{q\in S^c\setminus\{d\}} (I-U_{\epsilon_q})   
  \sum_{j_d= 1} ^\infty  \sum_{i_d= j_d}^\infty  U_{\epsilon_{d}}^{-j_d} g_{S,i_d} = \\
   \sum_{S\subset\{1,\dots,d\},\,d\in S} \prod_{q\in S^c} (I-U_{\epsilon_q})  g_{S} + 
   \sum_{S\subset\{1,\dots,d-1\}}  \prod_{q\in S^c} (I-U_{\epsilon_q})   U_{\epsilon_{d}}^{-j_d} g_{S} 
  \endmultline
$$
where for  $S\subset\{1,\dots,d-1\}$, 
$$
  g_S =  \sum_{j_d= 1}^\infty \sum_{i_d= j_d}^\infty U_{\epsilon_{d}}^{-j_d}  g_{S,i_d}\,\,\,\text{and }\,\,\,
   g_{S\cup \{d\}} = \sum_{i_d= 0}^\infty g_{S,i_d}.
$$
This proves  (2) with $g_S$ defined by (5b), for $d$ parameters; (4e) guarantees convergence.

\medskip

Now, let's suppose (2). For $d=1$, (4) follows by Proposition 3. Let us assume that the implication $(2) \implies (4)$  is true for $d-1$, $d\geq 2$. 
We consider the random field defined by $U_{\epsilon_1}, \dots, U_{\epsilon_{d-1}}$ and
$$
 f= \sum_{i_d=0}^\infty \Big(  \sum_{i_1= 0}^\infty \dots \sum_{i_{d-1}= 0}^\infty a_{i_1,\dots,i_{d-1},i_d}  U_{-i_1,\dots-,i_{d-1},0}
  (U_{\epsilon_d}^{-i_d}e)\Big);
$$
$i_d$ is a parameter. We thus have
$$
  f = \sum_{S\subset \{1,\dots,d\}}  \prod_{q\in S^c} (I-U_{\epsilon_q}) g_S = 
  \sum_{S\subset \{1,\dots,d-1\}}  \prod_{q\in S^c\setminus \{d\}} (I-U_{\epsilon_q}) \bar g_S
$$
where 
$$
  \bar g_S = g_{S\cup \{d\}} +  (I-U_{\epsilon_d}) g_S,\quad S\subset\{1,\dots,d-1\};
$$
by (5b) we have
$$
   \bar g_S = \sum_{i_d=0}^\infty  \sum_{i_r\geq 0, \,r\in S}  \sum_{j_u\geq 1,\,u\in S^c\setminus \{d\}} \sum_{i_u\geq j_u\,u\in S^c\setminus \{d\}}.  a_{i_r,i_u,i_d} U_{0,-j_u,-i_d} e
$$
By using Heyde's condition (6) we get, for $S\subset\{1,\dots,d-1\}$,
$$\multline
  \infty > \sum_{j_d\geq 1} \big\| \sum_{i_d\geq j_d} P_0^{(d)} U_{\epsilon_d}^{i_d} \bar g_S \big\|_2^2 = \\ 
  \sum_{j_d\geq 1} 
   \big\| \sum_{i_r\geq 0, \,r\in S}  \sum_{j_u\geq 1,\,u\in S^c\setminus \{d\}} \sum_{i_u\geq j_u\,u\in S^c\setminus \{d\}}  \sum_{i_d\geq j_d}   
  a_{i_r,i_u,i_d} U_{0,-j_u,0} e \big\|_2^2 =    \\
  \sum_{j_d\geq 1}  \sum_{j_u\geq 1,\,u\in S^c\setminus \{d\}} \big(  \sum_{i_r\geq 0, \,r\in S}  \sum_{i_u\geq j_u\,u\in S^c\setminus \{d\}} 
  \sum_{i_d\geq j_d} a_{i_r,i_u,i_d} \big)^2.
  \endmultline
$$
In (6) we can sum the $j_d$ from 0, hence we also have
$$
  \sum_{j_u\geq 1,\,u\in S^c\setminus \{d\}} \big(  \sum_{i_r\geq 0, \,r\in S\cup \{d\}}  \sum_{i_u\geq j_u\,u\in S^c\setminus \{d\}}  a_{i_r,i_u} \big)^2 <\infty;
$$
we thus have proved (4e). By the preceding implication we also get (5b).

In the same way we can prove the theorem for non adapted random fields.
The proof is thus accomplished.
\qed
\enddemo

\bigskip

In \cite{V06} it was proved that for $1\leq p<\infty$, convergence of (3)
in $L^p$ is a necessary and sufficient condition for the decomposition (1) in $L^p$. The result can be extended to $\Bbb Z^d$ actions.
Let us denote
$$
  Q_j^{(q)}f = E(f |  \Cal F_j^{(q)}) = \sum_{i\leq j} P_i^{(q)}f.
$$ 

\proclaim{Theorem 4} Let  $1\leq p<\infty$,  $f\in L^p$. (2) is equivalent to the convergence in $L^p$ of
$$\gathered
   \sum_{i_r\in \Bbb Z,\,r\in S} \sum_{i_u\geq 0,\,u\in S'}  \sum_{i_v\geq 1,\,v\in S^c\setminus S'}  
  \prod_{r\in S} \prod_{u\in S'} \prod_{v\in S^c\setminus S'}  P_0^{(r)} Q_{-1}^{(u)} (I- Q_{-1}^{(v)})
   U_{\epsilon_r}^{i_r}  U_{\epsilon_u}^{i_u} U_{\epsilon_v}^{-i_v}   f  = \\
   \sum_{i_r\in \Bbb Z,\,r\in S} \sum_{j_u\geq 1,\,i_u\geq 0,\,u\in S} \sum_{j_v\geq 0,\,i_v\geq 1,\,v\in S^c} P_{i_r,-j_u,j_v}U_{ i_u,-i_v} f
\endgathered \tag{3a}
$$
for all $S\subset \{1,\dots,d\}$. 

If $f$ is adapted then a necessary and sufficient condition for (2) is the  convergence in $L^p$ of
$$
  \sum_{i_1=0}^\infty \dots \sum_{i_d=0}^\infty E(U_{i_1,\dots,i_d}f \,|\, \Cal F_{\underline{0}}). \tag{3b}
$$
\endproclaim

Like in Theorem 2, (5) gives us the transfer functions $g_S$.

\demo{Proof} We prove the theorem for adapted functions only. For $f$ adapted, (3a) means the convergence of
$$
   \sum_{i_r\geq 0,\,r\in S} \sum_{i_u\geq 0,\,u\in S^c} \prod_{r\in S} \prod_{u\in S^c}  P_0^{(r)} Q_{-1}^{(u)}  
  U_{\epsilon_r}^{i_r}  U_{\epsilon_u}^{i_u} f,
$$
$S\subset \{1,\dots,d\}$. In the same way as in the proof of Theorem 2 we can prove equivalence with (3b).
A generalization to nonadapted (regular) random fields is straightforward.

Let us suppose the convergence in (3b). As shown in the Introduction, the functions $U_{\epsilon_1}^if$ are $\Cal F_0^{(q)}$-measurable for
$2\leq q\leq d$, hence $E(U_{\epsilon_1}^if | \Cal F_0^{(1)}) = E(U_{\epsilon_1}^if | \Cal F_{\underline{0}})$, $i\in \Bbb Z$.
By (3b) and Proposition 1 we thus get
$$
  f = m_1 + g_1- U_{\epsilon_1}g_1, \quad g_1 = \sum_{i=0}^\infty E(U_{\epsilon_1}^if | \Cal F_{\underline{0}}),\,\,\,
  m_1 = \sum_{i=0}^\infty P_0^{(1)} U_{\epsilon_1}^if;
$$
(3b) applies to $g_1, m_1$.
By iterating the procedure we get (2).
\medskip 

Now let us suppose that (2) is true. In the same way as in the proof of Theorem 2 we can see that $f = m_1 +  g_1- U_{\epsilon_1}g_1$
where $g_1, m_1$ can be decomposed by (2), and by Proposition 1 we have $g_1 = 
\sum_{i=0}^\infty E(U_{\epsilon_1}^if | \Cal F_{\underline{0}})$, $ m_1 = \sum_{i=0}^\infty P_0^{(1)} U_{\epsilon_1}^if$. After having repeated
the procedure for $g_1, m_1$ we will get the convergence of $\sum_{i_2=0}^\infty \sum_{i_1=0}^\infty 
E( U_{\epsilon_1}^{i_1}  U_{\epsilon_2}^{i_2} f | \Cal F_{\underline{0}})$ and by iterating the procedure we will get (3b).

\enddemo
\qed

In one dimensional case we can use $\sum_{i=0}^\infty E(f | \Cal F_{-1})$ or $\sum_{i=0}^\infty E(f | \Cal F_0)$; convergence of
the series is equivalent. The first series has the advantage of giving a cobounding function. In dimension $d\geq 2$ the convergence of
$$
  \sum_{i_1=0}^\infty \dots \sum_{i_d=0}^\infty E(U_{i_1,\dots,i_d}f \,|\, \Cal F_{-\underline{1}})
$$
is, however, not sufficient for (2). As an example we can consider a two dimensional adapted stationary random field (7) where $a_{i,j} =0$ for 
$i\geq 1$,  $a_{0,0}=0$, and $a_{0,j} = 1/j$ for $j\geq 1$.
\medskip

In what follows we will present two conditions which are sufficient (but not necessary) for (2).
The first  was proved by M\. El Machkouri and D\. Giraudo in  \cite{ElM16} and  is formulated for adapted functions. The result can, nevertheless,
be extended to the general (regular) case.

\proclaim{Theorem 5 (El Machkouri, Giraudo)} If $f$ is adapted and 
$$
  \sum_{i_1=0}^\infty \dots \sum_{i_d=0}^\infty \|E(U_{i_1,\dots,i_d}f \,|\, \Cal F_{\underline{0}})\|_p <\infty \tag{9}
$$
then the martingale-coboundary representation (2) holds.
\endproclaim

Theorem 5 extends Theorem 2 from \cite{Go69} to random fields. Its relation to Theorem 4 is the same as the relation of  
Theorem 2 from \cite{Go69} to Theorem 2 from \cite{V93}.

Let us denote 
$$
  \bar i = \cases i\,\,\,&\text{if}\,\,\, i\leq -1, \\
                          i+1 \,\,\,&\text{if}\,\,\, i\geq 0.
  \endcases
$$

\proclaim{Theorem 6} Let $p\geq 2$. If
$$
  \sum_{i_1\in\Bbb Z} \dots \sum_{i_d\in\Bbb Z}\bar  i_1^{2} \bar i_2^{2}\dots \bar i_d^{2}\|P_{i_1,i_2,\dots,i_d}f\|_p^2 <\infty \tag{10}
$$
then the martingale-coboundary decomposition (2) in $L^p$ holds.
\endproclaim

Before proving Theorem 6, let us prove

\proclaim{Lemma 7} Let $d\in \Bbb N$ and $a_{\underline{i}}$, $\underline{i} \in \Bbb N^d$,  be real numbers for which \newline
$ \sum_{i_1=1}^\infty\dots \sum_{i_d=1}^\infty  i_1^{2}i_2^{2}\dots i_d^{2} a_{\underline{i}}^2 < \infty$.
Then for every $\underline{j} \in \Bbb N^d$ the sum \newline
$\sum_{i_1=j_1}^\infty\dots \sum_{i_d=j_d}^\infty  a_{\underline{i}}$
converges and
$$
  \sum_{j_1=1}^\infty \dots \sum_{j_d=1}^\infty \big( \sum_{i_1=j_1}^\infty\dots \sum_{i_d=j_d}^\infty 
  a_{\underline{i}} \big)^2 \leq C  \sum_{i_1=1}^\infty\dots \sum_{i_d=1}^\infty  
   i_1^{2}i_2^{2}\dots i_d^{2} a_{\underline{i}}^2 \tag{11}
$$
where $C<\infty$ is a universal constant.
\endproclaim

\underbar{Remark 1.} From Lemma 7 we deduce that (10) implies $\sum_{\underline{i} \in \Bbb Z^d} \|P_{\underline{i}} f\|_p <\infty$, 
i.e\. Hannan's condition in $L^p$ (cf\. \cite{VWa14}). For $p=2$ and $d=1$ Hannan's condition is independent of decomposition (2) (cf\. \cite{V93}).
\medskip

\demo{Proof of Lemma 7} First, let us suppose $d=1$, $j\geq 1$, and $1/2 < \alpha <1$.
We then have
$$
  \big(\sum_{i= j}^\infty a_i\big)^2 = \big(\sum_{i= j}^\infty \frac1{i^\alpha} i^\alpha a_i \big)^2 \leq
  \big(\sum_{i=j}^\infty \frac1{i^{2\alpha}}\big) \big(\sum_{i= j}^\infty i^{2\alpha} a_i^2 \big)
  \leq \frac{c}{j^{2\alpha-1}} \sum_{i= j}^\infty i^{2\alpha} a_i^2
$$  
hence
$$
  \sum_{j\geq 1} \big(\sum_{i= j}^\infty a_i\big)^2 \leq 
  c \sum_{j=1}^\infty \frac1{j^{2\alpha-1}} \sum_{i=j}^\infty i^{2\alpha} a_i^2 
  = c \sum_{i=1}^\infty \big(\sum_{j=1}^i \frac1{j^{2\alpha-1}}\big) i^{2\alpha} a_i^2 \leq 
  C \sum_{i=1}^\infty i^{2} a_i^2
$$  
for some constants $c, C$.

For $d\geq 2$ we can prove (11) by induction. We present this just for $d=2$:
$$\gather
  \sum_{j_1=1}^\infty \sum_{j_2=1}^\infty \big( \sum_{i_1=j_1}^\infty \sum_{i_2=j_2}^\infty  a_{i_1,i_2} \big)^2
  \leq  \sum_{j_1=1}^\infty \sum_{j_2=1}^\infty \frac{c}{j_1^{2\alpha-1}} \sum_{i_1=j_1}^\infty i_1^{2\alpha}
  \big(\sum_{i_2=j_2}^\infty  a_{i_1,i_2} \big)^2 = \\
  c \sum_{j_2=1}^\infty \sum_{i_1=1}^\infty  i_1^{2\alpha} \big(\sum_{j_1=1}^{i_1} \frac{1}{j_1^{2\alpha-1}}\big)
  \big( \sum_{i_2=j_2}^\infty  a_{i_1,i_2} \big)^2 \leq 
  C \sum_{j_2=1}^\infty \sum_{i_1=1}^\infty  i_1^{2} \big( \sum_{i_2=j_2}^\infty  a_{i_1,i_2} \big)^2 = \\
  C \sum_{i_1=1}^\infty  i_1^{2} \sum_{j_2=1}^\infty \big( \sum_{i_2=j_2}^\infty  a_{i_1,i_2} \big)^2 \leq
  C' \sum_{i_1=1}^\infty \sum_{i_2=1}^\infty  i_1^{2}  i_2^{2}  a_{i_1,i_2}^2.
  \endgather
$$
\enddemo
\qed

\demo{Proof of Theorem 6} Let $d=1$, suppose that $f$ is adapted. Let $J\subset \Bbb Z_+ = \{0,1,\dots \}$. Recall
$$
  \sum_{j\in J} E(U^jf | \Cal F_{0}) = \sum_{i\leq 0} P_i \sum_{j\in J} E(U^jf | \Cal F_{0}) =
  \sum_{i\leq 0}  \big( \sum_{j\in J} P_i E(U^jf | \Cal F_{0}) \big) =
  \sum_{i\leq 0} \sum_{j\in J} P_i U^jf.
$$  
By Burkholder's inequality there is a $C_1$ such that
$$\gather
  \big\|\sum_{i\leq 0} \big(\sum_{j\in J} P_i U^jf \big)\big\|_p \leq
  C_1 \big\|\sum_{i\leq 0} \big(\sum_{j\in J} P_i U^jf \big)^2\big\|_{p/2}^{1/2} \leq 
  C_1\big( \sum_{i\leq 0} \big\|\sum_{j\in J} P_i U^jf\|_p^2\big)^{1/2} \leq \\
  C_1 \Big( \sum_{i\leq 0} \big( \sum_{j\in J} \| P_i U^jf\|_p \big)^2 \Big)^{1/2} = 
  C_1  \Big( \sum_{i=0}^\infty \big( \sum_{j\in i+J} \|P_0 U^jf\|_p \big)^2 \Big)^{1/2}
  \endgather
$$  
and by Lemma 7 there is a $C_2$ such that
$$
  \sum_{i=0}^\infty \big( \sum_{j=i}^\infty \|P_0 U^jf\|_p \big)^2 \leq  C_2 \sum_{i=0}^\infty \bar i^{2}  \|P_0 U^i f\|_p^2.
$$
From this we deduce the convergence of  $ \sum_{j= 0}^\infty E(U^jf | \Cal F_{0})$ (in $L^p$). 
In the same way we can prove the theorem for a (regular) non adapted $f$. 

By induction the result can be extended  to all $d\geq 1$. To show this, let us consider $d=2$.

Denote $f_j = P_{-j}^{(2)}f$, $j\geq 0$. Thus  $f=\sum_{j=0}^\infty f_j$ and  assumption (10) implies
$$
  \sum_{i=1}^\infty \sum_{j=1}^\infty  i^2j^2 \|P_{-i+1}^{(1)}f_{j-1} \|_p^2 <\infty. 
$$
In particular, for $j\geq 0$ we have $\sum_{i=1}^\infty i^2 \|P_{-i+1}^{(1)}f_j \|_p^2 <\infty$. From the proof of Theorem 6 for $d=1$
it follows that there exist $m_j, g_j \in L^p$, $j\geq 0$, such that 
$$\gather
  f_j = m_j + g_j- U_{\epsilon_1} g_j, \quad  P_{-j}^{(2)} m_j = m_j,\,\, P_{-j}^{(2)} g_j=g_j,\,\, P_0^{(1)} m_j =m_j; \\
  \|g_j\|_p^2 \leq C \sum_{i=1}^\infty i^2 \|P_{-i}^{(1)}f_j \|_p^2.
  \endgather
$$
Because $\|P_{-i}^{(1)}f_j \|_p = \|P_{0}^{(1)}U_{\epsilon_1}^if_j \|_p = \|P_{\underline{0}}U_{i,j}f\|_p$, and because by (10) we have
$ \infty > \sum_{i=1}^\infty \sum_{j=1}^\infty  i^2j^2 \|P_{\underline{0}}U_{i-1,j-1}f\|_p^2 $, we by Lemma 7 get 
$ \sum_{j=0}^\infty \|g_j\|_p <\infty$. Therefore the series $g= \sum_{j=0}^\infty g_j$ converges in $L^p$. 

Because $g_j =  P_{-j}^{(2)}g$ and $\sum_{j=1}^\infty j^2 \|g_{j-1}\|_p^2 <\infty$ we get a martingale-coboundary decomposition
$g = G - U_{\epsilon_2}G$, $G\in L^p$. 

In the same way we get $m = \sum_{j=0}^\infty m_j = M - U_{\epsilon_2} M$, $M\in L^p$.
This results in a martingale-coboundary representation (2) for $f$.

\enddemo
\qed

\underbar{Remark 2.} For  $p=2$, the assumption of Theorem 6 follows from Gordin - El Machkouri - Giraudo's condition (assumption of Theorem 5).

\demo{Proof}
It is sufficient to study the adapted case only because the proof of the general case is similar. Let
$f = \sum_{l=0}^\infty \sum_{i=0}^\infty a_{l,i} U^{-i}e_l$, $\sum_{l=0}^\infty \sum_{i=0}^\infty a_{l,i}^2 <\infty$, $P_0e_l=e_l$, $\|e_l\|_2=1$,
and $e_l$ are mutually orthogonal. 
Let us suppose
$$
  \sum_{j=0}^\infty \| E(U^jf | \Cal F_{0}) \|_2 = \sum_{j= 0}^\infty \big( \sum_{l=0}^\infty \sum_{i=j}^\infty a_{l,i}^2 \big)^{1/2} < \infty; 
$$
we have 
$$
  \sum_{i=0}^\infty \bar i^2 \|P_{-i}f\|_2^2 =   \sum_{l=0}^\infty \sum_{i= 0}^\infty \bar i^2  a_{l,i}^2
$$
where $\bar i = i+1$.
For $k\geq 0$ let us denote $b_k = \big(\sum_{l=0}^\infty  \sum_{i=k}^\infty a_{l,i}^2\big)^{1/2}$; we thus have $\| E(U^jf | \Cal F_{0}) \|_2 = b_j$.
Then 
$$
  \sum_{k=0}^\infty \bar kb_k^2 = \sum_{l=0}^\infty \sum_{i=0}^\infty \big(\sum_{k=1}^{\bar i} k\big) a_{l,i}^2 \geq 
  C \sum_{l=0}^\infty \sum_{i= 0}^\infty \bar i^2  a_{l,i}^2 \tag{12}
$$
and
$$
  \sum_{k=0}^\infty \bar k b_k^2 = \sum_{k=0}^\infty \bar k b_k^{3/2} b_k^{1/2} \leq 
  \sqrt{\sum_{k=0}^\infty \bar k^2 b_k^3} \sqrt{\sum_{k=0}^\infty b_k}.\tag{13}
$$
Let us suppose $\sum_{k=0}^\infty b_k <\infty$. By Abel's summation
$$
  \sum_{k=1}^n ((k+1)-k) b_k = (n+1)b_{n+1} - b_1 - \sum_{k=1}^n (k+1)(b_{k+1}-b_k). \tag{14}
$$
The sequence of $b_n$ is decreasing and the sums on the left converge. If $\sup_n nb_n = \infty$ then the series 
$\sum_{k=1}^n (k+1)(b_{k+1}-b_k)$ diverges and so $nb_n \to\infty$. This contradicts $\sum_{k=1}^\infty b_k <\infty$. Therefore, 
 $\sup_n nb_n < \infty$. From this and (13) it follows  $\sum_{k=1}^\infty k b_k^2 <\infty$.

Before passing to higher dimension let us notice that $\sum_{k=1}^\infty b_k <\infty$ implies $\liminf_{k\to\infty} kb_k =0$. In (14)
we thus get $\sum_{k=1}^\infty (k+1)(b_{k+1}-b_k) = b_1 + \sum_{k=1}^\infty b_k$ hence 
$$
  nb_n \leq 2b_1 + 2 \sum_{k=1}^\infty b_k \leq 4 \sum_{k=1}^\infty b_k, \quad n\geq 1. 
$$
\medskip

To present a proof for $d>1$, let us consider $d=2$ (for $d>2$ the proof is analogical).
We suppose that $f = \sum_{l=0}^\infty \sum_{i_1=0}^\infty  \sum_{i_2=0}^\infty  a_{l,i_1,i_2} U_{-i_1,-i_2}e_l$ where 
$e_l$ are mutually orthogonal,  $P_{0,0} e_l=e_l$, $\|e_l\|_2=1$, and $\sum_{l=0}^\infty  \sum_{i_1=0}^\infty  \sum_{i_2=0}^\infty  a_{l,i_1,i_2}^2 <\infty$.
Like in the one dimensional case we denote
$$
  b_{j_1,j_2}^2 =  \sum_{l=0}^\infty \sum_{i_1=j_1}^\infty  \sum_{i_2=j_2}^\infty a_{l,i_1,i_2}^2 = \| E(U_{j_1,j_2}f | \Cal F_{0,0})\|_2^2.
$$
We have 
$$
  \sum_{i_1=0}^\infty  \sum_{i_2=0}^\infty\bar  i_1^2 \bar i_2^2 \|P_{-i_1,-i_2}f\|_2^2 = 
  \sum_{l=0}^\infty \sum_{i_1=0}^\infty  \sum_{i_2=0}^\infty  \bar  i_1^2\bar  i_2^2 a_{l,i_1,i_2}^2
$$
and using the same idea as in (12) and (13) 
$$\multline
   \sum_{j_1=0}^\infty  \sum_{j_2=0}^\infty \bar j_1 \bar j_2  b_{j_1,j_2}^2 =  
  \sum_{l=0}^\infty \sum_{j_1=0}^\infty  \sum_{j_2=0}^\infty \bar j_1 \bar j_2
   \sum_{i_1=j_1}^\infty  \sum_{i_2=j_2}^\infty a_{l,i_1,i_2}^2 = \\
   \sum_{l=0}^\infty \sum_{j_1=0}^\infty  \sum_{i_1=j_1}^\infty \bar j_1 \big(\sum_{j_2=0}^\infty \bar j_2 \sum_{i_2=j_2}^\infty a_{l,i_1,i_2}^2 \big)
  \geq \sum_{l=0}^\infty \sum_{j_1=0}^\infty  \sum_{i_1=j_1}^\infty \bar j_1  \big(C \sum_{i_2=0}^\infty  \bar i_2^2  a_{l,i_1,i_2}^2 \big) \geq \\
  C^2 \sum_{l=0}^\infty \sum_{j_1=0}^\infty  \sum_{i_2=0}^\infty \bar  j_1^2 \bar j_2^2  a_{l,i_1,i_2}^2,
\endmultline
$$
$$\multline
   \sum_{j_1=1}^\infty  \sum_{j_2=1}^\infty \bar j_1 \bar j_2  b_{j_1,j_2}^2 =  
  \sum_{j_1=1}^\infty  \sum_{j_2=1}^\infty \bar j_1 \bar j_2 b_{j_1,j_2}^{3/2}b_{j_1,j_2}^{1/2} \leq \\
  \big( \sum_{j_1=1}^\infty  \sum_{j_2=1}^\infty \bar j_1^2 \bar j_2^2 b_{j_1,j_2}^3\big)^{1/2} 
  \big(  \sum_{j_1=1}^\infty  \sum_{j_2=1}^\infty  b_{j_1,j_2} \big)^{1/2}.
  \endmultline
$$
Let us suppose that $\sum_{j_1=0}^\infty  \sum_{j_2=0}^\infty b_{j_1,j_2} <\infty$.
In the same way as before we prove that for every $j_1\geq 1$, $j_1 \sup_{j_2\geq 1} j_2 b_{j_1,j_2} \leq 4j_1 \sum_{j_2=1}^\infty  b_{j_1,j_2}
= j_1 B_{j_1}$. Because $\sum_{j_1=1}^\infty B_{j_1} < \infty$ we in an analogical way as in the one dimensional case deduce $\sup_{j_1} j_1B_{j_1}<\infty$
hence $\sup_{j_1,j_2\geq 1} \bar j_1 \bar j_2  b_{j_1,j_2} <\infty$. Therefore, $\sum_{j_1=0}^\infty  \sum_{j_2=0}^\infty \bar j_1 \bar j_2  b_{j_1,j_2}^2  <\infty$.
   
\enddemo
\qed

\underbar{Remark 3.} There exists a dynamical system $(\Omega, \Cal A, T, \mu)$ with an increasing filtration $\Cal F_i = T^{-i}\Cal F$ 
and an $f\in L^2$ such that $f= \sum_{i=0}^\infty a_i U^{-i}e$, $\sum_{i=0}^\infty a_i^2 < \infty$, $e\in L^2$, and $(U^ie)$ is a martingale
difference sequence, 
$$
  \sum_{k=0}^\infty \|E(f | \Cal F_{-k})\|_2 = \infty, \,\,\, \text{and}\,\,\,
  \sum_{k=0}^\infty k^2 \|P_{-k}f\|_2^2 = \sum_{k=0}^\infty k^2 a_k^2 < \infty.
$$

\demo{Proof} Let $(U^ie)$ be a martingale difference sequence. 
 We define $n_k= 2^k$, $k\geq 0$, $\epsilon_0=1$, and $\epsilon_k = 2^{-k}/k$,  $k\geq 1$.
Let  $a_i= \epsilon_k>0$ for $i=n_k$ and $a_i=0$ for all other $i$. 
For $f= \sum_{i\geq 0} a_iU^{-i}e$ it then holds
$
  \sum_{i\leq 0} i^2 \|P_if\|_2^2 =  \sum_{i=0}^\infty i^2 a_i^2 = \sum_{k=0}^\infty n_k^2 \epsilon_k^2 <\infty
$
and 
$
  \sum_{i=0}^\infty \|E(f | \Cal F_{-i})\|_2 \geq \sum_{k=0}^\infty n_k \epsilon_k = \infty.
$
\enddemo
\qed

\bigskip
\subheading{3.Applications}
\bigskip

Martingale-coboundary decomposition has played an important role in the study of limit theorems for stationary sequences of random variables.
For random fields, the martingale-coboundary representation has proved useful as well.
\bigskip

\underbar{1. Invariance principle.} We are interested in the weak convergence of normalized sums $(1/|\underline{n}|) \sum_{\underline{1}\leq \underline{i}
\leq \underline{n.t}} U_{\underline{i}}f$, $\underline{n.t} = (n_1t,\dots,n_dt)$, to a Brownian sheet in $D[0, 1]^d$. We will call this  case 
the weak invariance principle (WIP).

For proving a WIP we first need a central limit theorem  for fields of martingale differences. 
If $d\geq 2$ the CLT, however, does not need to hold even in the case of an ergodic field of orthomartingale differences (cf\. \cite{WaWo13}).
The CLT is true if one of the generating transformations $T_{\epsilon_q}$ is ergodic, cf\. \cite{V15}. In particular, this assumption is valid if the 
$\Bbb Z^d$ action is
Bernoulli, i.e\. the $\sigma$-algebra $\Cal A$ is generated by iid random variables $e\circ T_{\underline{i}}$, $\underline{i}\in \Bbb Z^d$, and the
filtration is defined by $\Cal F_{\underline{i}} = \sigma\{e\circ  T_{\underline{j}}: \,\,\underline{j} \leq \underline{i}\}$. (The filtration
is completely commuting then.)

A WIP under Hannan's condition for random fields was proved in \cite{VWa14}.
In \cite{ElG16} El Machkouri and Giraudo proved that assumptions of Theorem 5 imply a WIP. The next theorem extends their result. 

\proclaim{Theorem 8} If the representation (2) applies and all components are square integrable then the WIP holds.
\endproclaim

The proof is, in fact, included in \cite{VWa14}, the proof of Theorem 5.1.  In that paper a  martingale-coboundary representation (2) with all terms in
$L^2$ is proved for $f=\sum_{\underline{i}} P_{\underline{i}}f$ where in the sum, only finitely many terms are nonzero (this follows from each of
the Theorems 2, 4, 5, 6 here). The WIP for such functions is proved and the proof uses square integrability and (2) only.


In \cite{V93} it is exposed that martingale-coboundary representation and Hannan's condition are independent. Therefore, Theorem 7 does not follow
from the WIP proved in \cite{VWa14}. 


For $d=1$ the martingale-coboundary representation is ``suboptimal" because WIP follows from Dedecker-Rio's and Maxwell-Woodroofe's 
conditions which are weaker. For $d\geq2$ we, nevertheless, do not know an analogue to the conditions of  Dedecker-Rio or Maxwell-Woodroofe. 
The version of Maxwell-Woodroofe's condition used in \cite{WaWo13} implies Hannan's condition.

\bigskip

\underbar{2. Probabilities of large deviations.} Let $(X_i)$ be a strictly stationary sequence of martingale differences, $X_i \in L^p$, $p\geq 2$. 
By \cite{D01, Proposition 1(a)} 
$$
  E\big|\sum_{i=1}^n X_i \big|^p \leq (2p)^{p/2} n^{p/2} \|X_1\|_p^p. \tag{15}
$$This improves  Burkholder inequality as presented in \cite{HaHe}
and also the inequality $ E\big|\sum_{i=1}^n X_i \big|^p \leq (18pq^{(1/2)})^p n^{p/2}  E|X_1|^p$ in \cite{LV01, p.150}.

Let $X_{\underline{i}}$, $\underline{1} \in \Bbb Z^d$, be strictly stationary orthomartingale differences. For $\underline{n} = (n_1,\dots,n_d)\in \Bbb N^d$,
$\sum_{i_1=1}^{n_1} \dots \sum_{i_{d-1}=1}^{n_{d-1}} X_{i_1,\dots,i_d}$, $i_d=1,\dots, n_d$, are stricly stationary martingale 
differences. Consequently,  by (15)
$$
  E \big| \sum_{i_d=1}^{n_d} \big( \sum_{i_1=1}^{n_1} \dots \sum_{i_{d-1}=1}^{n_{d-1}} X_{i_1,\dots,i_d}\big)\big|^p  \leq 
  (2p)^{p/2} n^{p/2} \big\|\sum_{i_1=1}^{n_1} \dots \sum_{i_{d-1}=1}^{n_{d-1}} X_{i_1,\dots,i_{d-1},1}\big\|_p^p.
$$
By iterating this approach we deduce
$$
  E \big|\sum_{i_1=1}^{n_1} \dots  \sum_{i_d=1}^{n_d}  X_{i_1,\dots,i_d} \big|^p  \leq 
   (2p)^{dp/2} \underline{n}^{p/2} \| X_{\underline{1}}\|_p^p. \tag{16}
$$

Let us suppose (2) with $g_S\in L^p$, $S\subset \{1,\dots,d\}$. Let  $S = \{q_1,\dots,q_r\} \subset \{1,\dots,d\}$.  Because  for each $q\in S$, $U_{\epsilon_q}^i g_S$ are martingale differences,
for $1\leq n_1, \dots, n_r < \infty$ it holds
$$
  E \big| \sum_{i_1=1}^{n_1}\dots \sum_{i_r=1}^{n_r} U_{\epsilon_{q_1}}^{i_r}\dots  U_{\epsilon_{q_r}}^{i_r} g_S \big|^p \leq
   (2p)^{rp/2} (n_1. \dots .n_r)^{p/2} \|g_S\|_p^p
$$ 
so that
$$
 \big\| \prod_{q\in S^c} (I-U_{\epsilon_q})g_S \big\|_p \leq (2p)^{r/2} 2^{d-r} |\underline{n}|^{1/2} \|g_S\|_p.
$$
We thus get the following theorem ( the second statement is deduced with the help of  Markov Inequality).

\proclaim{Theorem 9}
$$
 \big\| \sum_{\underline{1}\leq \underline{i} \leq \underline{n}} U_{\underline{i}}f \big\|_p \leq
  2^{2d}p^{d/2}  |\underline{n}|^{1/2} \sum_{S\subset \{1,\dots,d\}} \|g_S\|_p, \tag{17}
$$
$$
  \mu\big( \sum_{\underline{1}\leq \underline{i} \leq \underline{n}} U_{\underline{i}}f > x|\underline{n}| \big) \leq
   2^{2dp}p^{dp/2} \Big(\sum_{S\subset \{1,\dots,d\}} \|g_S\|_p \Big)^p x^{-p} |\underline{n}|^{-p/2}.  \tag{18}
$$ 
\endproclaim

For the dimension $d=1$ an estimate similar to (18) was found in  \cite{LV01}; stationarity is not needed there (only uniformly bounded
$L^p$ norms). Under the assumptions of Theorem 5, estimates similar to (17), (18) were proved in  \cite{ElG16}. \newline

Using the same ideas as in \cite{ElG16} the results can be extended to Orlicz spaces. 
Let us recall (cf\. \cite{KR}) that for a Young function $\psi$ (a real convex nondecreasing function on $\Bbb R_+$, $\psi(0)=0$, $\lim_{x\to\infty} \psi(x) =\infty$)
the Orlicz space $L_\psi$ associated to $\psi$ is the space of all random variables $Z$ such that for some $c>0$, $E(\psi(Z/c)) <\infty$. The 
Luxemburg norm $\|Z\|_\psi$ of $Z$ is then defined by $\|Z\|_\psi = \inf \{ c>0 : E(\psi(Z/c)) \leq 1\}$. $(L_\psi, \|.\|_\psi)$ is a Banach space.

Let us define $h_\alpha = ((1-\alpha)/\alpha)^{1/\alpha} 1_{\{0<\alpha<1\}}$ and a Young function 
$$
  \psi_\alpha(x) = \exp((x+h_\alpha)^\alpha) - \exp(h_\alpha^\alpha).
$$
In the same way as in \cite{ElVWu13, Lemma 4} we can deduce, for $ \psi_\alpha$ defined above, $0<q<2/d$, $\beta(q) = 2q/(2-dq)$, and for a positive 
constant $C$  depending only on $d$ and $q$, 
the following estimates.

\proclaim{Theorem 10}
$$\gather
   \big\| \sum_{\underline{1}\leq \underline{i} \leq \underline{n}} U_{\underline{i}}f \big\|_{\psi_q} \leq C | \underline{n}|^{1/2}
   \sum_{S\subset \{1,\dots,d\}} \|g_S\|_{\psi_{\beta_q}}, \\
  \mu\big( \sum_{\underline{0}\leq \underline{i} \leq \underline{n}} U_{\underline{i}}f > x \big) \leq
  (1+ e^{h_q^q}) \exp\big(- \big(\frac{x}{C |\underline{n}|^{1/2} \sum_{S\subset \{1,\dots,d\}} \|g_S\|_{\psi_{\beta_q}} +h_q}\big)^q\big).
\endgather
$$
\endproclaim

If $x=|\underline{n}|$ and $q=2/3$ we get an estimate of ordre $ C \exp(-c |\underline{n}|^{1/3})$.
For $d=1$ it is $\mu(S_n \geq n) \leq C \exp(-c n^{1/3})$ which was found, like (18), in \cite{LV01} where only uniform boundednes 
of moments is needed.
In \cite{LV01} it is proved  that  for strictly stationary sequences of martingale differences both estimations are (up to a constant and for $d=1$) optimal.

\bigskip

More applications can be found in the field of  reverse martingale approximation for noninvertible commuting transformations, cf\. \cite{DeGo14}
and \cite{ CDV15}.



\enddocument

\end